\documentclass[12pt]{amsart}

\usepackage{amsmath,amssymb}
\usepackage{amsfonts}
\usepackage{amsthm}
\usepackage{latexsym}
\usepackage{graphicx}
\usepackage{epsfig}

\newtheorem{lemma}{Lemma}

\newtheorem{theorem}{Theorem}

\newenvironment{remark}[1][Remark]
           {\medbreak\noindent {\em #1. \enspace}}
           {\par \medbreak}

\makeatletter \@addtoreset{equation}{section} \makeatother

\def\ddt{\frac{d}{dt}}

\begin{document}

\title{Cross Curvature Flow on Locally Homogenous Three-manifolds (I) }

\author{Xiaodong Cao$^*$}
\thanks{$^*$Research
partially supported by an MSRI postdoctoral fellowship}
% by NSF
%grant no. \# }

\address{Department of Mathematics,
  Cornell University, Ithaca, NY 14853}
\email{cao@math.cornell.edu}
%\curraddr{MSRI}
%\thanks{}

\author{Yilong Ni}
%\thanks{$^2$Research
%partially supported by NSF grant no. \# }

\address{Department of Mathematics,
  University of Oklahoma, Norman, OK 73019}
\email{yni@aftermath.math.ou.edu}
%\curraddr{}
%\thanks{}

%    Information for third author
\author{Laurent Saloff-Coste$^{\flat}$}
\thanks{$^{\flat}$Research
partially supported by NSF grant no. \# DMS 0603886}

\address{Department of Mathematics,
  Cornell University, Ithaca, NY 14853}
\email{lsc@math.cornell.edu}
%\curraddr{}
%\thanks{}

%    General info
\renewcommand{\subjclassname}{%
  \textup{2000} Mathematics Subject Classification}
% Global differential geometry
% 53C21 Methods of Riemannian geometry, including PDE methods;
%   curvature restrictions [See also 58J60]
% 53C44 Geometric evolution equations (mean curvature flow)
% 53C55 Hermitian and K\"ahlerian manifolds [See also 32Cxx]
% Qualitative properties of solutions
% 35B35 Stability, boundedness
%Parabolic equations and systems [See also 35Bxx, 35Dxx, 35R30, 35R35, 58J35]
%35K55 Nonlinear PDE of parabolic type
%35K90 Abstract parabolic evolution equations
% Partial differential equations on manifolds; differential operators
%58J35 Heat and other parabolic equation methods
% 58J37 Perturbations; asymptotics
% Systems theory; control - Stability
% 93D05 Lyapunov and other classical stabilities
%       (Lagrange, Poisson, $L^p, l^p$, etc.)

\date{Nov. 6th,  2007}

\maketitle

\markboth{Xiaodong Cao, Yilong Ni and Laurent Saloff-Coste} {Cross
Curvature Flow }

\begin{abstract} Chow and Hamilton introduced the
cross curvature flow on closed $3$-manifolds with negative or
positive sectional curvature. In this paper, we study the negative
cross curvature flow in the case of  locally homogenous metrics on
$3$-manifolds. In each case, we describe the long time behavior of
the solutions of the corresponding ODE system.
\end{abstract}

\section{Introduction}

\subsection{Homogeneous evolution equations}

Hamilton's Ricci flow (\cite{H3}) is the seminal and most
successful example of the idea of deforming a Riemannian structure
using a geometric evolution equation. Special cases arise when the
metric is invariant under a group of transformations and this
property is preserved by the flow. In particular, if the group of
isometries of the original Riemannian structure is transitive,
then the geometric evolution equation reduces to an ODE in the
tangent space of an arbitrary  fixed origin. In this spirit, the
Ricci flow on locally homogeneous $3$-manifolds was analyzed by
Isenberg and Jackson [1992]
%\cite{isenberg1}
, quasi-convergence of model geometries under the Ricci flow was
studied by Knopf and McLeod [2001]
%\cite{knopf2}
, and the case of the Ricci flow on locally homogenous closed
$4$-manifolds was analyzed by Isenberg, Jackson and Lu [2006]
%\cite{isenberg2}
. Lott [2007]
%\cite{Lott}
interprets these results using
the notion of groupoids and solitons.

This paper studies the asymptotic behavior of the (negative) cross
curvature flow on locally homogeneous metrics in dimension $3$.
This flow was introduced by Chow and Hamilton [2004]
%\cite{chowhamiltonxcf}
and is (so far) specific to dimension $3$.
It depends on a sign choice (see Section 1.3 below). Chow and
Hamilton conjectured that for any compact $3$-manifold that admits
a metric with negative sectional curvature, the (positive) cross
curvature flow exists for all time and converge to a hyperbolic
metric. Because of the structure of cross curvature flow equation,
no general  existence results are expected when sectional
curvatures take different signs, which is the case for most
homogeneous geometries. However, in the case of homogeneous
geometries, there is no difficulties in defining, say, the
negative  curvature flow since the equations reduce to a system of
ODEs. The choice of a sign mentioned above can then be interpreted
simply as running the flow either forward of backward (although
one should observe that it is not clear, a priori, which direction
should be considered the forward direction).

These remarks prompted us to study the asymptotic forward {\em
and} backward behaviors  of the maximal two sided solutions of the
cross curvature flow on homogeneous $3$-manifolds. The present
paper deals with the (forward) behavior of the negative cross
curvature flow. The companion paper \cite{csc2} deals with the
backward behavior of the negative cross curvature flow (i.e., the
forward behavior of  the positive  cross curvature flow). The
backward behavior of the Ricci flow will be considered elsewhere.
One interesting discovery is that, generically, the backward
behavior of these flows is described by convergence (of the
distance functions) to non-degenerate sub-Riemannian geometries.
See \cite{csc2}. Concerning the forward direction studied in this
paper, it is interesting to compare the behavior of the Ricci flow
to that of the negative cross curvature flow. Let us  briefly
describe the similarities and differences. In the case of
geometries modeled on  $\mbox{SU}(2)$ and the Heisenberg group,
the behavior of the Ricci flow and cross curvature flow are
similar. On $\mbox{SU}(2)$, both  flows (in their normalized
version) have solutions that exists for all positive time and
converge towards  round metrics. On closed manifolds with
Heisenberg type geometries both normalized flows exist for all
positive time and as $t$ tends to infinity, they produce almost
flat metrics. On closed manifolds of type $E(2)$, both normalized
flows exists for all positive time. The normalized  Ricci flow
converges to a flat metric whereas the normalized negative cross
curvature flow produces almost flat metrics but develops a cigar
degeneracy. On $E(1,1)$ type closed manifolds (i.e., Sol
geometries), the normalized Ricci flow exists for all positive
time and presents a cigar degeneracy whereas the normalized cross
curvature flow exists only for a finite time and there is a
dimensional collapse with the sectional curvatures blowing up. The
case of compact quotients of $SL(2,\mathbb R)$ is the most
difficult and perhaps the most interesting. The normalized Ricci
flow exists for all time and presents a pancake type degeneracy.
For the normalized negative cross curvature flow, two different
types of behavior are possible. For  metrics with a specific
symmetry, the flow exists for all time and develops a pancake
degeneracy. For generic (homogeneous) metrics, the flow exists for
a finite time, there is a dimensional collapse and the curvatures
blow up (to plus and minus infinity).

In the rest of this introduction, we quickly review the necessary material
on locally homogenous $3$-manifolds as well as the
definition of the cross curvature flow. Sections 2 to 6 are devoted to
the different examples: Heisenberg, $E(1,1)$ (i.e., Sol), $\mbox{SU}(2)$,
$SL(2,\mathbb R)$ and $E(2)$.

\subsection{The cross curvature tensor on $3$-manifolds}

On a $3$-dimensional Riemannian manifold $(M,g)$, let $Rc$ be the
Ricci tensor and $R$ be the scalar curvature. The Einstein tensor
is defined by $E=Rc-\frac{1}{2}Rg$, and its local components are
$E_{ij}=R_{ij}-\frac{1}{2}Rg_{ij}.$  Raising the indices, define
$P^{ij}=g^{ik}g^{jl}R_{kl}-\frac{1}{2}Rg^{ij},$ where $g^{ij}$ is the
inverse of $g_{ij}$. Let $V_{ij}$ be the inverse of $P^{ij}$ (if
it exists). The cross curvature tensor is
$$h_{ij}=\left(\frac{\det P^{kl}}{\det g^{kl}}\right)V_{ij}.$$
This definition is taken from \cite{chowhamiltonxcf}.

Assume that computation are done in an orthonormal frame where the
Ricci tensor is diagonal. Then the cross curvature tensor is
diagonal.  If the principal sectional curvatures are $k_1,k_2,k_3$
($k_i=K_{jkjk}$, circularly) so that
$R_{ii}=k_j+k_l \;\;(\mbox{circularly in $i,j,l$}),$
then
\begin{equation} \label{C=kk}
h_{ii} = k_jk_l.
\end{equation}
Notice that this definition works even when some of the sectional
curvatures vanish (this was also addressed in \cite{ma06}).

The following lemma is proved in \cite{chowhamiltonxcf} using the
contracted second Bianchi identity,

\begin{lemma} (Chow and Hamilton [2004]
%\cite{chowhamiltonxcf}
). The cross curvature tensor $h$ satisfies the following
identities,
\[
\nabla_i P^{ij}=0 \;\;\mbox{ and }\;\;
(h^{-1})^{ij} \nabla_i h_{jk}=\frac{1}{2} \nabla_k h_{ij}.
\]
Moreover, if the sectional curvatures are either positive or
negative, then the identity map
$id: (M,h)\rightarrow (M,g)$
is a harmonic map.
\end{lemma}

Recall that when the Ricci curvature tensor is positive
(negative), the identity map
$id: (M,g)\rightarrow (M,\pm Rc)$
is a harmonic map.

\subsection{The cross curvature flows}
%In \cite{chowhamiltonxcf} ,
Chow and Hamilton [2004] define the cross curvature flow on
$3$-manifolds having either positive sectional curvature or
negative sectional curvature.  The local existence of the flow,
under any one of these two circumstances, was proved by Buckland
[2006]
%\cite{buckland}
. More precisely, if $\epsilon =\pm 1$ is the
sectional curvature sign of the metric $g_0$, the cross curvature
flow starting from $g_0$ is the solution of
$$\left\{\begin{array}{l}\frac{\partial}{\partial
t}g=-2\epsilon h,\\
g(0)=g_0.\end{array}\right.$$

For the purpose of this paper, it should be noticed that locally
homogeneous manifolds seldom have sectional curvatures that are
all of the same sign. In dimension 3, positive sectional curvature
is only possible on locally homogeneous manifolds covered by the
sphere $\mbox{SU}(2)$. Negative sectional curvature occurs only on
hyperbolic $3$-manifolds. All other locally homogeneous closed
Riemannian $3$-manifold are either flat or have some positive
sectional curvatures \cite[Theorem 1.6]{milnor76}. Thus the
definition above is not really practical for our purpose. In fact,
at least in the case of locally homogeneous $3$-manifolds, it
seems very natural to investigate both the positive and the
negative cross curvature flows where the positive cross curvature
flow is defined by
$$\left\{\begin{array}{l}
\frac{\partial}{\partial t}g=2 h,\\
g(0)=g_0.\end{array}\right.\leqno{(+\mbox{XCF})}$$ and the
negative cross curvature flow is defined by
$$\left\{\begin{array}{l}
\frac{\partial}{\partial t}g= -2 h,\\
g(0)=g_0.\end{array}\right. \leqno{(-\mbox{XCF})}$$
In this paper we consider the negative cross curvature flow ($-$XCF).

As in the Ricci flow, we can also consider the normalized cross
curvature flow (NXCF) on closed $3$-manifolds. It preserves the
volume of closed 3-manifolds and is given by
%$$
%\frac{\partial}{\partial
%t}g_{ij}=2h_{ij}-\frac{2}{3}\bar{h}g_{ij}, \leqno{(+\mbox{NXCF})}
%$$
$$
\frac{\partial}{\partial
t}g_{ij}=-2h_{ij}+\frac{2}{3}\bar{h}g_{ij}, \leqno{(\mbox{NXCF})}
$$
where $\bar{h}=\int_{M^3}g^{ij}h_{ij}du/\int_{M^3}du$. As for the
Ricci flow, the flows (-XCF) and (NXCF) only differ by a change of
scale in space and a re-parametrization of time
($\tilde{g}(\tilde{t})=\psi (t)g(t)$, $\tilde{t}=\int \psi^2(t)$).

\subsection{Locally homogeneous  $3$-manifolds}
Following Isenberg and Jackson [1992]
%\cite{isenberg1}
(to which we refer for details concerning the following
discussion), we take the view point that our original interest is
in closed Riemannian $3$-manifolds that are locally homogeneous.
By a result of Singer [1960]
%\cite{singer60}
, the universal cover of a locally homogeneous
manifold is homogeneous, that is, its isometry group acts
transitively. Now, since the cross curvature flow (just as the
Ricci flow) commutes with the projection map from the universal
cover, we can as well study the flow on the (often non-compact)
universal cover.

In dimension $3$ there are 9 possibilities for the universal
cover, which can be labelled by the minimal isometry group that
acts transitively:
\begin{itemize}
\item[(a)]$H(3)$ ($H(n)$ denotes the isometry group of hyperbolic
$n$-space); $\mbox{SO}(3)\times \mathbb R$; $H(2)\times \mathbb
R$; \item[(b)] $\mathbb R^3$; $\mbox{SU}(2)$; $\mbox{SL}(2,\mathbb
R)$; Heisenberg; $E(1,1)=\mbox{Sol}$ (the group of isometry of
plane with a flat Lorentz metric); $E(2)$ (the group of isometries
of the Euclidean plane).
\end{itemize}
The crucial difference between cases (a) and (b) above is that, in
case (b), the universal cover of the corresponding closed
$3$-manifold is (essentially) the minimal transitive group of isometries itself
(with the caveat that both $\mbox{SL}(2,\mathbb R)$ and $E(2)$
should be replaced by their universal cover)
whereas in case (a) this minimal group is
of higher dimension. The cases (a) and (b) are studied separately.
The case (b) is called the Bianchi case in \cite{isenberg1}. It
corresponds exactly to the classification of $3$-dimensional
simply connected unimodular Lie groups (nonunimodular Lie groups
cannot cover a closed manifold).

\subsection{Real $3$-dimensional unimodular Lie algebra}\label{Milnor-f}
The basis for our study is Milnor's description in \cite[Section
4]{milnor76} of all three dimensional real Lie algebras equipped
with an Euclidean structure (i.e., a left-invariant metric $g_0$
on the Lie group). Remember that the data here is the Lie algebra with
a fixed Euclidean structure (and, in fact, a fixed orientation).
The crucial result is as follows. Assume that $\mathfrak g$ is a
$3$-dimensional real unimodular Lie algebra equipped with an
oriented Euclidean structure. Then there exists a (positively
oriented) orthonormal basis $(e_1,e_2,e_3)$ and reals
$\lambda_1,\lambda_2,\lambda_3$ such that the bracket operation of
the Lie algebra has the form
$$[e_i,e_j]=\lambda_k e_k  \;\;\;\mbox{(circularly in $i,j,k$)}.$$
Milnor shows that such a basis diagonalizes the Ricci tensor and
thus also the cross curvature tensor. If $f_i= a_ja_k e_i$ with
nonzero $a_i,a_j,a_k\in \mathbb R$, then $[f_i,f_j]= \lambda_k
a_k^2 f_k$ (circularly in $i,j,k$). Using the choice of
orientation, we may assume that at most one of the $\lambda_i$ is
negative and then, the Lie algebra structure is entirely
determined by the signs (in $\{-1,0,+1\}$) of
$\lambda_1,\lambda_2,\lambda_3$ as follows:
$$\begin{array}{cccl}
+&+&+&  \;\mbox{SU}(2)\\
+&+&-& \; \mbox{SL}(2,\mathbb R) \\
+&+&0& \;  E(2) \; \mbox{(Euclidean
motions in $2D$)} \\
+&0&-&  \;E(1,1) \; \mbox{(also called Sol)}\\
+&0&0& \;  \mbox{Heisenberg group} \\
0&0&0& \;  \mathbb R^3
\end{array}$$
In each case, let $\epsilon=(\epsilon_1,\epsilon_2,\epsilon_3)\in
\{-1,0,+1\}^3$ be the corresponding choice of signs. Then, given $\epsilon$
and an Euclidean  metric $g_0$ on the corresponding Lie algebra,
we can choose a basis $f_1,f_2,f_3$ (with $f_i$
collinear to $e_i$ above) such that
\begin{equation}\label{MF}
[f_i,f_j]= 2\epsilon_k f_k \;\;\;\mbox{(circularly in $i,j,k$)}.
\end{equation}
As mentioned above, the metric, the Ricci tensor and the cross curvature tensor
are diagonalized in this
basis and this property is obviously maintained throughout either
the Ricci flow or cross curvature flow. We call $(f_i)_1^3$ a
Milnor frame for $g_0$. If we let $(f^i)_1^3$ be the dual frame
of $(f_i)_1^3$, the metric $g_0$ is diagonalized by this frame and
has the form
\begin{equation}\label{g_0}
g_0= A_0 f^1\otimes f^1 +B_0 f^2\otimes f^2+
C_0 f^3\otimes f^3.\end{equation} Assuming existence of the
flow $g(t)$ starting from $g_0$, under either the Ricci flow or
the cross curvature flow (positive or negative), the original
frame $(f_i)_1^3$ stays a Milnor frame for $g(t)$ along the flow.
Thus, $g(t)$ has the form
\begin{equation}\label{gflow}
g(t)= A(t) f^1\otimes f^1 +B(t)f^2\otimes f^2+ C(t) f^3\otimes
f^3.\end{equation} It follows that these flows reduce to  ODEs in
$(A,B,C)$. Given a flow, the explicit form of the ODE depends on
the underlying Lie algebra structure. With the help of the
curvature computations done by Milnor in [1976]
%\cite{milnor76}
, one can find the explicit form of the equations
for each Lie algebra structure. The Ricci flow case was treated in
\cite{isenberg1}. The case of the negative cross curvature flow is
treated below.

\subsection{The trivial cases}
The three non-Bianchi cases and the flat case $\mathbb R^3$ all
lead to essentially trivial behaviors. For $\mathbb R^3$, this is
obvious.

In the hyperbolic case $H(3)$, the only homogeneous metrics
are constant multiple of the standard hyperbolic metric.
They all have constant negative curvature. The cross curvature
tensor is a multiple of the identity. So each metric is a fixed
point under the NXCF in this case.

In the case of $\mbox{SO}(3)\times \mathbb R$, the homogeneous
metrics must have a product form corresponding to  a metric on
$\mathbb R$ and a round
metric on the $2$ sphere. In a proper frame, two of the principal
sectional curvatures vanish and thus $h=0$. The cross curvature
flow is trivial.

Finally, for $H(2)\times \mathbb R$, the metrics again have a
product form so that two of the three sectional curvatures vanish
and $h=0$. The cross curvature flow is trivial.

\subsection{ Acknowledgement} The authors would like to thank
Professor Bennett Chow and David Glickenstein for their interest
and helpful conversations concerning the cross curvature flow.

\section{The negative XCF on the Heisenberg group (Nil geometries)}

Given a left-invariant  metric $g_0$ on the Heisenberg group,  fix a Milnor
frame $\{f_i\}_1^3$ such that
$$[f_2,f_3]=2f_1, \;\;[f_3,f_1]=0,\;\;[f_1,f_2]=0$$
and (\ref{g_0})-(\ref{gflow}) hold.  Using
\cite{milnor76}, the sectional curvatures are:
$$
K(f_2 \wedge f_3)= -\frac{3A}{BC},\;\; K(f_3 \wedge
f_1)=\frac{A}{BC},\;\; K(f_1 \wedge f_2)= \frac{A}{BC}.$$ The
scalar curvature is $R = -2A/(BC)$. The computation of the cross
curvature tensor easily follows by (\ref{C=kk}). In the frame
$(f_i)_1^3$ and its dual frame $(f^i)_1^3$, the cross curvature
tensor is given by
$$h=
  \frac{A^3}{B^2C^2}  f^1\otimes f^1
   -3\frac{A^2}{BC^2} f^2\otimes f^2
    -3\frac{A^2}{B^2C} f^3\otimes f^3 .$$
Hence, the negative cross
curvature flow ($-$XCF) reduces to the ODE system
$$\frac{dA}{dt} = -\frac{2A^3}{B^2C^2},\; \;\;
\frac{dB}{dt} =\frac{6A^2}{BC^2},\;\;\;
\frac{dC}{dt} = \frac{6A^2}{B^2C} .$$
Observe that
\[
\frac{dA}{Adt} = -2\frac{A^2}{B^2C^2} = -\frac{1}{3}\frac{dB}{Bdt}
= -\frac{1}{3}\frac{dC}{Cdt}.
\]
Hence  $B/C$, $A^3B$ and $A^3C$ stay constant under the flow
and
\[
\frac{dA}{dt} = -2\frac{A^3}{B^2C^2} =- 2A^3 \cdot
\frac{A^6}{A_0^6 B_0^2} \cdot \frac{A^6}{A_0^6 C_0^2}
=-\frac{2}{A_0^{12}B_0^2C_0^2} A^{15}.
\]
As $R_0=-2A_0/(B_0C_0)$, we arrive at $
A(t)   = A_0
(1+ 7R_0^2t)^{-\frac{1}{14}}$,
$B(t) = B_0
(1+7R_0^2t)^{\frac{3}{14}}$, and
$C(t) = C_0 (1+7 R_0^2t)^{\frac{3}{14}}$.
This shows that the solution of the flow exists for all time $t\ge 0$.
The sectional curvatures are
\[
K(f_2 \wedge f_3) =
\frac{3R_0}{2}(1+7R_0^2t)^{-\frac{1}{2}} \mbox{ and }
K(f_1 \wedge f_2) =K(f_3 \wedge f_1) =
-\frac{R_0}{2}(1+7R_0^2t)^{-\frac{1}{2}}.
\]

Hence we have the following result.
\begin{theorem}
On the Heisenberg group, for any initial data $A_0$, $B_0$, $C_0
>0$, the solution of the negative (XCF) on $[0,\infty)$ is given
by
$$A(t)   = A_0 (1+
7R_0^2t)^{-\frac{1}{14}},\;\; B(t) = B_0 (1+7R_0^2t)^{\frac{3}{14}}\;
\mbox{ and }\;C(t) = C_0 (1+7 R_0^2t)^{\frac{3}{14}},$$ where
$R_0=-2A_0/(B_0C_0)$. The sectional curvatures decay as
$t^{-1/2}$.
\end{theorem}

A closed manifold is a Nilmanifold if it is the quotient of a
nilpotent Lie group by a discrete subgroup. A closed Riemannian
manifold $(M,g)$ is $\epsilon$-flat ($\epsilon>0$ fixed) if it
admits a metric such that all sectional curvatures are bounded
above in absolute value by $\epsilon d^{-2}$ where $d$ is the
diameter of $(M,g)$. A manifold is almost flat if it admits
$\epsilon$-flat metrics for all small $\epsilon>0$. By a Theorem
of M. Gromov [1978]
%\cite{gromov}
(see also \cite{buserkarcher}),
in any dimension, a manifolds is almost flat if and only
if it is covered by a Nilmanifold.

The closed locally homogeneous $3$-manifolds associated to the
Heisenberg group are Nilmanifolds and thus are almost flat. Let
$d(t)$ be the diameter of such a manifold under the negative XCF
$g(t)$ considered above. Obviously, $d(t)^2$ is of order
$t^{3/14}$ and the sectional curvatures are bounded in absolute
value by  a constant times $t^{-1/2}$. This shows that, as $t$
tends to infinity, the negative XCF yields $\epsilon$-flat metric
(with $\epsilon(t)$ of order $t^{-2/7}$). The  normalized flow (NXCF), has
a similar behavior with a slightly different numerology.

\section{The negative XCF on Sol geometry (E(1,1))}

Given a left-invariant metric $g_0$ on $E(1,1)$,
fix a Milnor frame $\{f_i\}_1^3$ such
that
$$[f_2,f_3]=2f_1, \;\;[f_3,f_1]=0,\;\;[f_1,f_2]=-2f_3.$$

The sectional curvatures are:
\[
K(f_2 \wedge f_3) = \frac{(A-C)^2-4A^2}{ABC},
\]

\[
K(f_3 \wedge f_1) = \frac{(A+C)^2}{ABC},
\]

\[
K(f_1 \wedge f_2) = \frac{(A-C)^2-4C^2}{ABC}.
\]

In the frame $\{f_1, f_2, f_3\}$, we have

\[
(h_{ij}) = \left(%
\begin{array}{ccc}
  -\frac{A(A+C)^3(3C-A)}{(ABC)^2}&  &  \\
  &\frac{B(3A-C)(3C-A)(A+C)^2}{(ABC)^2} &  \\
   &  & -\frac{C(A+C)^3(3A-C)}{(ABC)^2} \\
\end{array}%
\right),
\]
 and the negative cross curvature flow equations are

\begin{equation}
\left \{
\begin{aligned}
\frac{dA}{dt} =& 2\frac{A(A+C)^3(3C-A)}{(ABC)^2},
\\
\frac{dB}{dt} =& -2\frac{B(3A-C)(3C-A)(A+C)^2}{(ABC)^2},
\\
\frac{dC}{dt} =& 2\frac{C(A+C)^3(3A-C)}{(ABC)^2}.
\end{aligned}
\right .
\end{equation}

If $A=C$ at $t=0$, then $A(t) \equiv C(t) $ as long as the solution exists.
Moreover,
\[
\frac{dB}{dt} = -2\frac{4A^2 \cdot 4A^2}{A^2A^2B} = -\frac{32}{B},
\]
so $B^2=B_0^2-64t$, that is, $B=\sqrt{B_0^2-64t}$.  Also, we have

\[
\frac{d \ln A}{dt} = \frac{32}{B_0^2-64t}.
\]
Hence
\[
A(t)=C(t) = \frac{A_0B_0}{ \sqrt{B_0^2-64t}}.
\]

For the case $A_0\ne C_0$, we may assume without loss of
generality that $A_0 > C_0$. Then we immediately have that $C$ is
increasing. Observing that
\begin{align*}
&\frac{d(A-C)}{dt}=-2\frac{(A+C)^4}{(ABC)^2}(A-C),\\
&\frac{d\ln(A/C)}{dt}=-8\frac{(A+C)^3}{(ABC)^2}(A-C),\\
&\frac{d(A-3C)}{dt}=-2\frac{(A+C)^3}{(ABC)^2}(A^2+6AC-3C^2),
\end{align*}
we find that $A>C$ and $A-C$, $A/C$ and $A-3C$ are decreasing.

Let us further assume that $3C_0>A_0$. Then we have $1<A/C<A_0/C_0<3$
and
$$
B\frac{dB}{dt}=-2\frac{(A+C)^2}{(AC)^2}(3A-C)(3C-A)
\in(-128,-E_0),
$$
where $E_0:=16(3C_0-A_0)/A_0>0$. Therefore there exists
$T_0\in(0,\infty)$ such that $B(T_0)=0$. Furthermore, when
$t\in[0,T_0)$
$
2E_0(T_0-t)<B^2<256(T_0-t)$.
Hence
$$
\frac1C\frac{dC}{dt}=2\frac{(A+C)^3}{(ABC)^2}(3A-C)>\frac{16}{B^2}
>\frac1{16(T_0-t)},
$$
which implies that $C$ and $A$ go to $\infty$ as $t\to T_0^-$.
It follows that,  as $t\to T_0^-$,
$$
B\sim \sqrt{64(T_0-t)},\quad A, C \sim
\frac{E_1}{\sqrt{T_0-t}},\quad A-C \sim E_2\sqrt{T_0-t},
$$
where $E_1, E_2$ are positive constants.

If $3C_0\le A_0$ then we claim that there
exists $t_1\ge0$, such that $3C(t_1)>A(t_1)$ (and thus $3C(t)>A(t)$
for all $t>t_1$, as long as the solution exists since $A-3C$ is decreasing).
Suppose on
the contrary that $3C\le A$ as long as the solution exists, then
we have that $B$ is increasing, $A$ is decreasing and
\begin{align*}
B\frac{dB}{dt}&=2\frac{(A+C)^2}{(AC)^2}(3A-C)(A-3C)<6\frac{(A+C)^2}{C^2}\\
&=6\left(1+A/C\right)^2<6\left(1+A_0/C_0\right)^2:=E_3.
\end{align*}
Therefore the solution exists for all $t\in[0,\infty)$ and $B^2<2E_3t+B_0^2$.
Furthermore,
\begin{align*}
\frac{d(A-3C)}{dt}&=-2\frac{(A+C)^3}{(ABC)^2}(A^2+6AC-3C^2)<
-\frac{16C^3\cdot 4C^2}{(ABC)^2}\\
&=-\frac{64C^3}{A^2B^2}<-\frac{64C_0^3}{A_0^2}\cdot\frac{1}{2E_3t+B_0^2}.
\end{align*}
Integrating the above inequality from $0$ to $\infty$ yields a contradiction.
Hence, we have the following theorem.
\begin{theorem}
On $E(1,1)$, for any initial data $A_0$, $B_0$, $C_0>0$, there
exists a time $T_0>0$, such that the solution of the negative cross
curvature flow exists for all $0\le t<T_0$.
Moreover, as $t\to T_0^-$,
$$
B\sim \sqrt{64(T_0-t)},\quad A, C \sim
\frac{E_1}{\sqrt{T_0-t}},\quad A-C \sim E_2\sqrt{T_0-t},
$$
where $E_1$ and $E_2$ are constants. The sectional curvatures
approach to $\pm\infty$ at rate of $(T_0-t)^{-1/2}$ as $t
\rightarrow T_0$.
\end{theorem}

\begin{remark}
Under the normalized flow, the solution also only exists up to a finite time $T_1$,
and
$$
B\sim E'_1(T_1-t),\quad A, C \sim \frac{E'_2}{\sqrt{T_1-t}}.
$$
The sectional curvatures approach to $\pm \infty$ at rate of
$(T_1-t)^{-1/2}$ as $t \rightarrow T_1$.
The diameter $d(t)$ increase like
$(T_1-t)^{-1/4}$, so the absolute values of the sectional
curvature are not $o(d(t)^{-2})$ (compare with the case of the Nil geometry).
Recall that the solution to the normalized Ricci flow exists for
all time and approaches a cigar degeneracy (see \cite{isenberg1}),
i.e., two directions shrink to zero while the other one expands to
$\infty$, and the sectional curvatures decay at rate of $t^{-1}$.
The Ricci flow and cross curvature flows behave quite differently
in this case.
\end{remark}

\section{The negative XCF on $SU(2)$}
Given a left-invariant metric $g_0$ on $SU(2)$, fix a Milnor frame
such that $[f_i,f_j]=2f_k$ circularly. We have $K(f_2 \wedge
f_3)=\frac{(B-C)^2}{ABC}-\frac{3A}{BC}+\frac2B+\frac2C$ and the
other sectional curvatures are obtained by circular permutation.
The cross curvature tensor is diagonal under the associated
orthogonal frame $\{f_i\}_1^3$  with $h_{11} = (ABC)^{-2}AYZ $ and
the other entries obtained again by circular permutation with
\begin{align*}
X=&3A^2-(B-C)^2-2AB-2AC, \\
Y=&3B^2-(A-C)^2-2AB-2BC, \\
Z=&3C^2-(A-B)^2-2BC-2AC.
\end{align*}
Therefore, under ($-$XCF), $A,B,C$ satisfy the following
equations
\begin{equation}\label{pdesu2}
\left \{
\begin{aligned}
\frac{dA}{dt}=&-2\frac{AYZ}{(ABC)^2},\\
\frac{dB}{dt}=&-2\frac{BZX}{(ABC)^2},\\
\frac{dC}{dt}=&-2\frac{CXY}{(ABC)^2}.
\end{aligned}
\right .
\end{equation}
Without loss of generality we may assume that $A_0\ge B_0\ge C_0$. Then we know
that $A\ge B\ge C$ as long as a solution exists. Observing that
\begin{align*}
Y&=(B-A)(A+B+2B-2C)-C^2\le -C^2<0,\\
Z&=-(A-B)^2+3C^2-2AC-2BC\le -C^2<0,
\end{align*}
we have
\begin{align*}
&\frac{d(A-B)}{dt}=\frac{2Z}{(ABC)^2}(A-B)(A^2+A(6B-2C)+(B-C)^2)\le 0,\\
&\frac{d(A-C)}{dt}=\frac{2Y}{(ABC)^2}(A-C)((A-B)^2+6AC-2BC+C^2)\le 0,\\
&\frac{d\ln(A/B)}{dt}=\frac{8Z}{(ABC)^2}(A-B)(A+B-C)\le 0,\\
&\frac{d\ln(A/C)}{dt}=\frac{8Y}{(ABC)^2}(A-C)(A+C-B)\le 0.
\end{align*}
It follows that $A$, $A-B$, $A-C$, $\ln(A/B)$ and $\ln(A/C)$ are decreasing.
Furthermore,
$$
-\frac{dA}{dt}=\frac{2AYZ}{(ABC)^2}\ge
\frac{2AC^4}{(A^2C)^2}=\frac{2C^2}{A^3} \ge
\frac{C_0^2}{A_0^2}\frac2A,
$$
which implies $\ddt A^2\leq -4C_0^2A_0^{-2}$. Therefore there
exists $T_0\in(0,\infty)$, such that $A(T_0)=B(T_0)=C(T_0)=0$. On
the other hand
$$
-\frac{dA}{dt}=\frac{2AYZ}{(ABC)^2}\le
\frac{2A}{(AC^2)^2}(3A^2)(4A^2) =\frac{24A^3}{C^4}\le
\frac{24A_0^4}{C_0^4}\frac1A.
$$
It follows that on $[0,T_0)$,
\begin{equation}\label{abound}
\sqrt{48(T_0-t)}\frac{A_0^2}{C_0^2}\ge A(t)\ge
2\sqrt{T_0-t}\frac{C_0}{A_0}.
\end{equation}
Since $A/C$ is decreasing and bounded below by $1$, we may assume that
$\lim_{T_0^-}A/C=p.$
We claim that $p=1$. Suppose instead that $p>1$. Then we have
$$
-\frac{d\ln(A/C)}{dt}=\frac{-4Y}{(ABC)^2}(A-C)(A+C-B)\ge
\frac{4C^2(A-C)C}{(A^2C)^2}\ge (1-p^{-1})\frac{C_0}{A_0}\frac4{A^2}.
$$
Integrating from $0$ to $T_0$ and using (\ref{abound}), we get a contradiction.
Therefore,
$\lim_{T_0^-}A/C=1. $
It follows easily from (\ref{pdesu2}) that as $t\to T_0^-$,
$
A,B,C\sim 2\sqrt{T_0-t}.
$
\begin{theorem}
For any choice of initial data $A_0$, $B_0$, $C_0>0$, there exists
a time $T_0>0$, such that the solution of the cross curvature flow
on $SU(2)$ exists for all $0\le t<T_0$. Moreover, as $t\to T_0^-$,
$$
A,B,C\sim 2\sqrt{T_0-t}.
$$
\end{theorem}

\begin{remark}
If we consider the normalized XCF, we find that solutions of
(NXCF) exists for all time. As $t\to \infty$, $A$, $B$ and $C$
approach a constant
 and the sectional curvatures also become constant, so the metric
becomes round. Since the sectional curvatures become positive for
$t$ large enough, according to Chow and Hamilton, the negative XCF
is the more natural choice in this case.
\end{remark}

%%%%%%%%%%%%%%%%%%%%%%%%%%%%%%%%%%%%%%%%%%%%%%%%%%%%%%%%%%%%%%%%%%%%%%%%%%%%%%%%%%%%%%%%%%%%%%%%%%%%%%%%%%%%%%%%%%%%%%%%%%%%%%%%

\section{The negative XCF on $SL(2,{\mathbb R})$}

Given a left-invariant metric $g_0$ on $SL(2,{\mathbb R})$, fix a
Milnor frame $\{f_i\}_1^3$ such that
$$[f_2,f_3]=-2f_1, \;\;[f_3,f_1]=2f_2,\;\;[f_1,f_2]=2f_3.$$
 The sectional curvatures are
\begin{align*}
K(f_2 \wedge f_3)&=\frac{1}{ABC}(-3A^2+B^2+C^2-2BC-2AC-2AB),\\
K(f_3 \wedge f_1)&=\frac{1}{ABC}(-3B^2+A^2+C^2+2BC+2AC-2AB),\\
K(f_1 \wedge f_2)&=\frac{1}{ABC}(-3C^2+A^2+B^2+2BC-2AC+2AB),
\end{align*}

and the cross curvature tensor under the associate
frame $(f_i)_1^3$ is
$$
(h_{ij}) = \frac1{(ABC)^2}\left(
\begin{array}{ccc}
  AF_2F_3 &  &  \\
  & B F_3F_1 &  \\
   &  & C F_1F_2 \\
\end{array}
\right),
$$

where
\begin{align*}
F_1=&-3A^2+B^2+C^2-2BC-2AC-2AB, \\
F_2=&-3B^2+A^2+C^2+2BC+2AC-2AB, \\
F_3=&-3C^2+A^2+B^2+2BC-2AC+2AB.
\end{align*}

Therefore, under the negative XCF, $A,B,C$ satisfy the following
equations
\begin{equation}\label{pdesl2}
\left \{
\begin{aligned}
\frac{dA}{dt}=&-\frac{2AF_2F_3}{(ABC)^2},\\
\frac{dB}{dt}=&-\frac{2BF_3F_1}{(ABC)^2},\\
\frac{dC}{dt}=&-\frac{2CF_1F_2}{(ABC)^2}.
\end{aligned}
\right .
\end{equation}
If $B_0=C_0$, then $B=C$ as long as a solution exists and $A, B$ satisfy
$$
\frac{dA}{dt}=-\frac{2A^3}{B^4},\quad
\frac{dB}{dt}=2\frac{3A^2+4AB}{B^{3}}.
$$
Then $A$ is decreasing and $B$ is increasing. We have
\begin{align*}
\frac{dA^{-2}}{dt}&=-2A^{-3}\frac{dA}{dt}\le 4B_0^{-4},\\
\frac{dB^3}{dt}&\le6\frac{3A_0^2+4A_0B}{B}\le C_1,
\end{align*}
where $C_1=24A_0+18A_0^2B_0^{-1}$ is a constant. Integrating from
$0$ to $t$, we obtain that
\begin{align*}
A&\ge(4B_0^{-4}t+A_0^{-2})^{-\frac12},\\
B&\le(C_1t+B_0^3)^{\frac13}.
\end{align*}
It follows that a solution exists on $[0,\infty)$.
\begin{align}\label{beqc}
&\frac{d(4A^{-1}+B^{-1})}{dt}=-4A^{-2}\frac{dA}{dt}-B^{-2}\frac{dB}{dt}
=-6\frac{A^2}{B^5}, \nonumber\\
&\frac{d(A^9B^3)}{dt}=9A^8B^3\frac{dA}{dt}+3A^9B^2\frac{dB}{dt}=24A^{10}.
\end{align}
Hence $4A^{-1}+B^{-1}$ is decreasing, which implies
$\lim_{t\to\infty}A:=A_{\infty}>0$. Integrating (\ref{beqc}) we
obtain that $A^9B^3\to \infty$ as $t\to\infty$. Hence
$\lim_{t\to\infty} B=\infty$. It is not hard to show that as
$t\to\infty$,
$$
A\sim A_\infty+\frac{1}{8 \sqrt[3]{3}}
A_\infty^{\frac53}t^{-\frac13}, \qquad B\sim (24A_\infty
t)^{\frac13}.
$$
For the case $B_0\ne C_0$, we may assume without loss of
generality that $B_0>C_0$. Then $B>C$ as long as a solution
exists. It follows that
$$
F_3=(B-C)(2A+B+3C)+A^2>A^2>0.
$$
Let $a=AB^{-1}$ and $c=CB^{-1}$.
\begin{lemma}
Suppose that at $t=0$, $a$ and $c$ satisfy
\begin{equation}\label{f2n}
a<1-c+2\sqrt{1-c}
\end{equation}
and
\begin{equation}\label{f1n}
a>\frac13(2\sqrt{1-c+c^2}-1-c).
\end{equation}
Then $a$ and $c$ satisfy (\ref{f2n}) and (\ref{f1n}) as long
as a solution exists.
\end{lemma}
\begin{proof} As
\begin{align*}
&F_2=(A-(B-C+2\sqrt{(B-C)B}))(A-(B-C-2\sqrt{(B-C)B})),\\
&F_1=(B-(A+C+2\sqrt{(A+C)A}))(B-(A+C-2\sqrt{(A+C)A})).
\end{align*}
we see that (\ref{f1n}) is equivalent to $F_1<0$ and
(\ref{f2n}) is equivalent to $F_2<0$. Since
$$
\left.\frac{dA}{dt}\right|_{F_2=0}=0,\quad
\left.\frac{dB}{dt}\right|_{F_2=0}>0\quad\mbox{ and }
\left.\frac{dC}{dt}\right|_{F_2=0}=0
$$
we obtain that
\begin{equation}\label{df2dt}
\left.\frac{dF_2}{dt}\right|_{F_2=0}<0.
\end{equation}
To prove that $F_2(t)<0$ we argue by contradiction. Suppose $t_0$
is the first time such that $F_2(t_0)=0$. Since $F_2(0)<0$, we
know that $F_2'(t_0)\ge 0$, which contradicts (\ref{df2dt}).
Therefore $F_2(t)<0$, which is equivalent to (\ref{f2n}).
Similarly $\left.\frac{dF_1}{dt}\right|_{F_1=0}<0$ and $F_1(t)<0$.
This completes the proof of the lemma.
\end{proof}
\begin{lemma}
Suppose that at $t=0$, $a$ and $c$ satisfy
\begin{equation}\label{f2p}
a\ge 1-c+2\sqrt{1-c}
\end{equation}
Then eventually $a$ and $c$ will satisfy (\ref{f2n}) and (\ref{f1n}).
\end{lemma}
\begin{proof}
Suppose on the contrary that $a\ge 1-c+2\sqrt{1-c}$ always holds, then
$F_2>0$ and thus $F_1<0$. Hence
$A$ is decreasing and $B,C$ are increasing as long as a solution exists.
Note that
$$
\frac{d\ln(C/B)}{dt}=\frac{8(B-C)}{(ABC)^2}(A+B+C)F_1<0.
$$
Therefore $c=C/B$ is decreasing, which implies that
$$
A=aB\ge(1-c+2\sqrt{1-c})B\ge (1-c_0+2\sqrt{1-c_0})B_0.
$$
On the other hand
$$
\frac{dB^2}{dt}=\frac{-4F_1}{A^2}\cdot\frac{F_3}{C^2}.
$$
As
$$F_3/C^2\ge 3((B/C)-1)\ge 3((B_0/C_0)-1)$$
and
$$-F_1=F_2+2((A+B)^2-C^2) \ge B^2-C^2 \ge C_0^2((B_0/C_0)^2-1),$$
it follows that $(B^2)'>\eta$ for some positive constant $\eta$.
If $B$ stays finite for all $t$ then the solution
exists on $[0,\infty)$ and $\lim_{t\to\infty}B=\infty$. Therefore as
$t\to\infty$
$$
F_2=-(B-C)(2A+C+3B)+A^2\le -3BC(\frac BC-1)+A^2 \to -\infty.
$$
This contradicts $ F_2>0$. If $B$ goes to $\infty$ in finite time then $F_2\to
-\infty$ in finite time. We also get a contradiction.
\end{proof}
\begin{lemma}
Suppose that at $t=0$, $a$ and $c$ satisfy
\begin{equation}\label{f1p}
a\le \frac13(2\sqrt{1-c+c^2}-1-c).
\end{equation}
Then eventually $a$ and $c$ will satisfy (\ref{f2n}) and (\ref{f1n}).
\end{lemma}
\begin{proof}
Suppose on the contrary that $a\le\frac13(2\sqrt{1-c+c^2}-1-c)$ always holds.
Then $F_1>0$ and thus  $F_2<0$. It follows that
$A,C$ are increasing and $B$ is decreasing as long as a solution exists.
Since
$$
-F_2=F_1+2( (A+B)^2-C^2 )\ge 2A^2+4AB \mbox{ and } F_3\ge A^2,
$$
we obtain that $A'\ge 4A^3(BC)^{-2}\ge 4A_0^3B_0^{-4}$. If $A$
stays finite for all $t$ then the solution exists on $[0,\infty)$
and $\lim_{t\to\infty}A=\infty$. Therefore as $t\to\infty$
$$
F_1=(B-C)^2-A(3A+2C+2B)\le (B_0-C_0)^2-3A^2\to -\infty.
$$
This is a contradiction. If $A$ goes to $\infty$ in finite time then $F_1\to
-\infty$ in finite time. We also get a contradiction.
\end{proof}
From the above three lemmas, we can assume without loss of generality that
(\ref{f2n}) and (\ref{f1n}) hold at $t=0$. Then $F_1<0$, $F_2<0$, $A,B$ are
increasing and $C$ is decreasing as long as a solution exists.

\medskip
\begin{lemma}\label{lemma5}
Suppose that in addition we have $A_0+C_0\le B_0$, then there exists
$T>0$ such that
\begin{equation}\label{sim}
A,B\sim E(T-t)^{-\frac{1}{2}}, \quad C\sim 8\sqrt{T-t}\qquad\mbox{
as }t\to T^{-},
\end{equation}
where $E$ is a positive constant.
\end{lemma}
\begin{proof}
We first claim that $A+C\le B$ holds for all $t$. In fact we have
\begin{equation}\label{dlna}
\frac{d\ln(A/B)}{dt}=-\frac{8(A+B)}{(ABC)^2}F_3(A+C-B),
\end{equation}
which implies $\left.\frac{d(A/B)}{dt}\right|_{A+C-B=0}=0$. On the
other hand $\left.\frac{d(C/B)}{dt}\right|_{A+C-B=0}<0$. It
follows that
$$
\left.\frac{d(A+C-B)}{dt}\right|_{A+C-B=0}<0,
$$
which implies $A+C\le B$ as long as a solution exists. Therefore
$a=AB^{-1}$ is increasing and $-F_2\ge(A+B)^2-C^2\ge 2A(A+B+C)$. It follows
that
$$
\frac{d\ln A}{dt}=\frac{2(-F_2)F_3}{(ABC)^2}> \frac{4AB\cdot
A^2}{(ABC)^2} =\frac{4A}{BC^2}\ge \frac{4A_0}{B_0C_0^2}.
$$
 If the solution exists on $[0,\infty)$, integrating the
above inequality from $0$ to $\infty$ yields
$\lim_{t\to\infty}A=\infty$. Since $A+C\le B$, we also obtain that
$\lim_{t\to\infty}B=\infty$. Let $p:=\lim_{t\to\infty}A/B$. Then
$p$ must be $1$, otherwise integrating (\ref{dlna}) from $0$ to
$\infty$ we get a contradiction. Since $\lim_{t\to\infty}A/B=1$,
as $t\to\infty$, $F_1,F_2\sim -(A+B)^2$. It follows from
$$
\frac{dC^2}{dt}=-4\frac{F_1F_2}{A^2B^2},
$$
that $C$ goes to $0$ in finite time. This contradicts the
assumption that a solution exists on $[0,\infty)$. Therefore any
solution blows up in finite time. Suppose $[0,T)$ is the maximal
time interval of a solution. Let $C_T =\lim_{t\to T^{-}}C$. If
$C_T>0$, then it follows easily from (\ref{pdesl2}) that $A,B$
stay bounded on $[0,T)$. Therefore we may extend the solution
beyond $T$, which contradicts the maximality of the
time interval $[0,T)$. Hence $\lim_{t\to T^{-}}C=0$. From
$$
\frac{d(AC)}{dt}=-F_2\frac{4AC}{(ABC)^2}(B^2-(A+C)^2)\ge 0,
$$
we obtain that $A\to\infty$ and $B\ge A+C\to\infty$ as $t\to T^{-}$. Let
$p:=\lim_{t\to T^{-}}AB^{-1}$. Since $A+C\le B$, we have $p\le 1$. If $p<1$,
then as $t\to T^{-}$,
$$
\frac{d(A-B)}{dt}=\frac{2F_3}{(ABC)^2}(BF_1-AF_2)
=\frac{2F_3}{(ABC)^2}((A+B)^2(B-A-2C)+(B-A)C^2)>0,
$$
which is impossible since $A-B\sim -(1-p)B$. Therefore
$\lim_{t\to T^{-}}AB^{-1}=1$ and
$$
\lim_{t\to T^{-}}A^{-2}F_1=\lim_{t\to T^{-}}A^{-2}F_2=-4
\mbox{ and }\lim_{t\to T^{-}}A^{-2}F_3=4.
$$
Then (\ref{sim}) follows easily from (\ref{pdesl2}).
\end{proof}
The only case left now is that $A,B,C$ satisfy $A_0+C_0>B_0$ and
$F_1(0), F_2(0)<0$. If there is a time $T^*$, such that
$A(T^*)+C(T^*)=B(T^*)$, then from Lemma \ref{lemma5} we know the
behavior of the solution. If on the other hand $A+C>B$ for all
$t$, then
$$
-F_1\ge (A+B)^2-C^2
$$
and (\ref{dlna}) implies that $AB^{-1}$
is decreasing. It follows that
$$
-F_2=(B-C+2\sqrt{B(B-C)}-A)(A-(B-C-2\sqrt{B(B-C)}))\ge C_1 AB,
$$
where $C_1$ is some constant depending only on $A_0,B_0$ and $C_0$. Using the
above two inequalities and (\ref{pdesl2}) we obtain that
$$
\frac{dC^2}{dt}=-4\frac{F_1F_2}{(AB)^2}\le
-4\frac{C_1((A+B)^2-C^2)}{AB} \le -8C_1.
$$
Therefore, there exist a finite $T>0$ such that $C(T)=0$. Form
$$
\frac{d(BC)}{dt}=-F_1\frac{4BC}{(ABC)^2}(A^2-(B-C)^2)\ge 0,
$$
we obtain that $B\to \infty$ and $A>B-C\to\infty$ as $t\to T^{-}$. Again let
$p:=\lim_{t\to T^{-}}AB^{-1}$. Since $A+C>B$, we have $p\ge 1$. If $p>1$,
then as $t\to T^{-}$,
$$
\frac{d(A-B)}{dt}=\frac{2F_3}{(ABC)^2}(BF_1-AF_2)
=\frac{2F_3}{(ABC)^2}((A+B)^2(B-A-2C)+(B-A)C^2)<0,
$$
which is impossible since $A-B\sim (p-1)B$. Therefore
$\lim_{t\to T^{-}}AB^{-1}=1$,
$$
\lim_{t\to T^{-}}A^{-2}F_1=\lim_{t\to T^{-}}A^{-2}F_2=-4
\mbox{ and }\lim_{t\to T^{-}}A^{-2}F_3=4
$$
and we have (\ref{sim}). Hence we have the following theorem.
\begin{theorem}
On $SL(2,R)$, for given initial data $A_0$, $B_0$, $C_0>0$, if
$B_0=C_0$, then $B(t)=C(t)$, and the solution of the negative
cross curvature flow exists for all $t\in[0,\infty)$. As functions
of $t$, $A$ is decreasing and $\lim_{t\to\infty}A=A_\infty>0$,
whereas  $B,C$ are increasing and go to $\infty$. Moreover, as
$t\to\infty$
$$
A\sim A_\infty+\frac{1}{8 \sqrt[3]{3}}
A_\infty^{\frac53}t^{-\frac13}, \qquad B=C \sim (24A_\infty
t)^{\frac13},
$$ the sectional curvatures all approach to $0$ as $t\rightarrow
\infty$ ($K(f_2 \wedge f_3) \sim -\frac{4}{B} \sim -E_1 t^{-1/3}$
and $K(f_3 \wedge f_1)=K(f_1 \wedge f_2)=\frac{A}{BC} \sim E_2
t^{-2/3}$, where $E_1$ and $E_2$ are some constants).

 If $B_0>C_0$, then there exists a time $T_0>0$, such that
the solution of the cross curvature flow on $SL(2,{\bf R})$ exists
for all $0\le t<T_0$. Moreover, as $t\to T_0^-$
$$
A,B\sim E(T_0-t)^{-\frac{1}{2}}, \quad C\sim 8\sqrt{T_0-t},
$$
where $E$ is some constant, and all sectional curvatures go to
$\pm\infty$ at the rate of $(T_0-t)^{-1/2}$ as $t\rightarrow T_0$.
\end{theorem}

\begin{remark}
The asymptotic behavior of the solution depends on the initial
data in this case, i.e., the condition $B=C$ does not characterize
the typical geometry of general solution. If we consider the
normalized cross curvature flow, then for the case of $B=C$,
$A\sim E_1 t^{-\frac25}$, $B=C \sim E_2 t^{\frac15}$, where $E_1$
and $E_2$ are some constants, one sectional curvature decays at
rate of $t^{-\frac{1}{5}}$ and the other two sectional curvatures
decay at rate of $t^{-\frac{4}{5}}$, we have a pancake degeneracy.
For the case of $B\neq C$ under (NXCF), $ A,B\sim E_1
(T_1-t)^{-\frac{1}{4}}$, $C\sim E_2 \sqrt{T_1-t}$, all sectional
curvatures go to $\pm\infty$ at the rate of $(T_1-t)^{-1/2}$ as
$t\rightarrow T_1$, where $T_1$ is the maximal existence time for
the solution. Recall that the solution of the Ricci flow in this
case exists for all time and develops a pancake degeneracy.
\end{remark}

\section{The negative XCF on $E(2)$}
Given a left-invariant metric $g_0$, fix a Milnor frame $\{f_i\}_1^3$ such
that
$$[f_2,f_3]=2f_1, \;\;[f_3,f_1]=2f_2,\;\;[f_1,f_2]=0.$$
The sectional curvatures are
\[
K(f_2 \wedge f_3) = \frac{1}{ABC} (B-A)(B+3A),
\]

\[
K(f_3 \wedge f_1) = \frac{1}{ABC} (A-B)(A+3B),
\]

\[
K(f_1 \wedge f_2) = \frac{1}{ABC} (A-B)^2,
\]
and the cross curvature tensor in the
frame $\{f_i\}_1^3$ is
$$
(h_{ij}) = \frac1{(ABC)^2}\left(
\begin{array}{ccc}
  AYZ &  &  \\
  & BZX &  \\
   &  & CXY \\
\end{array}
\right),
$$
where
$$X=(A-B)(3A+B),\;
Y=(B-A)(3B+A) \mbox{ and }
Z=-(A-B)^2.$$

Therefore, under the negative XCF, $A,B,C$ satisfy the following
equations
%\begin{equation}\label{pdee2}
%\left\{
%\begin{aligned}
%\frac{dA}{dt}=&-\frac{2AYZ}{(ABC)^2},\\
%\frac{dB}{dt}=&-\frac{2BZX}{(ABC)^2},\\
%\frac{dC}{dt}=&-\frac{2CXY}{(ABC)^2}.
%\end{aligned}
%\right .
%\end{equation}
%So
\begin{equation}\label{pdee2}
\left\{ \begin{aligned}
\frac{dA}{dt}=&-\frac{2A(3B+A)(A-B)^3}{(ABC)^2}, \\
\frac{dB}{dt}=&-\frac{2B(3A+B)(B-A)^3}{(ABC)^2}, \\
\frac{dC}{dt}=&\frac{2C(3A+B)(3B+A)(A-B)^2}{(ABC)^2}.
\end{aligned}\right.
\end{equation}

If $A_0=B_0$, then the geometry stays flat at all time. Without loss
of generality, we assume that $A_0>B_0$. Then we have $B_0\leq
B(t)<A(t)\leq A_0$ as long as the solution exists.
Since $C'(t)>0$, $C(t)$ is increasing.
It follows easily from (\ref{pdee2}) that a solution exists for all
$t\in[0,\infty)$. We first claim that $\lim_{t\to\infty}C=\infty$. In fact,
suppose that $\lim_{t\to\infty}C=C_\infty<\infty$, then
$$
\frac{d(A-B)}{dt}=-2\frac{(A-B)^3}{(ABC)^2}(A^2+6AB+B^2)
$$
implies $((A-B)^{-2})'\sim E_1$ as $t\to\infty$, where $E_1$ is
some constant. It follows that $A-B\sim
E_1^{-\frac12}t^{-\frac12}$. Then, from
$$
\frac{dC^2}{dt}=4\frac{(3A+B)(3B+A)}{(AB)^2}(A-B)^2,
$$
we obtain that $C\to\infty$ as $t\to\infty$. This is a
contradiction. Now,
$$
\frac{d \ln((A-B)^2C)}{dt}=2\frac{d \ln(A-B)}{dt}+\frac{d \ln
C}{dt} =\frac{2(A-B)^4(A+B)}{(ABC)^2}\le E_2\frac{dB}{dt},
$$
for some positive constant $E_2$. Therefore $(A-B)^2C$ is increasing and
approaches some finite number as $t\to\infty$. Hence
$$
C^2\frac{dC}{dt}=2\frac{(3A+B)(3B+A)}{(AB)^2}(A-B)^2C\sim E_3,
$$
as $t\to\infty$. It follows that $C\sim E_4t^{\frac13}$,
$A-B\sim E_5t^{-\frac16}$ and
$$
\frac{d(A+B)}{dt}=-2\frac{(A-B)^4(A+B)}{(ABC)^2}\sim E_6 t^{-\frac43}.
$$
Therefore we have the following theorem.
\begin{theorem}
On $E(2)$, for any initial data $A_0$, $B_0$, $C_0>0$, if
$A_0=B_0$, then the solution of ($-$XCF) exists for all
time, $A(t)=B(t)=A_0$ and $C(t)=C_0$ for all time $t$ (the
geometry stays flat).

If $A_0>B_0$, then the solution exists for all $t\in[0,\infty)$ and,
as $t\to \infty$
$$
A\sim E_1+E_2t^{-\frac16}, B\sim E_1-E_2t^{-\frac16} \mbox{ and }
C\sim (8E_2/E_1)\sqrt{6}t^{\frac13},
$$
where $E_1$ and $E_2$ are positive constants. Two of the sectional
curvatures decay like $t^{-1/2}$, while the other one decays like
$t^{-2/3}$.
\end{theorem}

\begin{remark}
Under the Ricci flow, the geometry converges to a flat metric. If
we consider the solution to the normalized cross curvature flow,
in the case that $A_0=B_0$, we still have flat metric. For the
case of $A_0\neq B_0$, we have
$$
A\sim E_1 t^{-\frac17}, B\sim E_1 t^{-\frac17} \mbox{ and } C\sim
E_2 t^{\frac27},
$$
while two of the sectional curvatures decay like $t^{-1/2}$, and
the other one decays like $t^{-5/7}$, hence the solution of
(NXCF) develops a cigar degeneracy, i.e., two directions shrink to
zero, the other one expands without bound, while the sectional
curvature dies off.
\end{remark}

\bibliographystyle{halpha}

%\bibliography{bio}
\end{document}